\documentclass[11pt]{article}
\usepackage{amsfonts}
\usepackage{amsmath,amssymb}
\begin{document}

\begin{center}
\large\bf Some addition to the generalized Riemann-Hilbert problem
\end{center}
\medskip
\begin{center}
R.R. Gontsov\footnote{Institute for Information Transmission
Problems, Moscow, Russia, rgontsov@inbox.ru.} and I.V.
Vyugin\footnote{Moscow State University, Russia,
ilya\_vyugin@mail.ru}
\end{center}

\begin{abstract}
We consider the generalized Riemann-Hilbert problem for linear
differential equations with irregular singularities. After
recalling the formulation of the problem in terms of vector
bundles over the Riemann sphere, we give some estimates for the
unique non-minimal Poincar\'e rank of the system and the number of
apparent singularities of the scalar equation constructed by
corresponding generalized monodromy data.
\end{abstract}
\medskip

\centerline{\bf\S1. Introduction}
\medskip

Consider a system
\begin{eqnarray}\label{syst}
\frac{dy}{dz}=B(z)y, \qquad y(z)\in{\mathbb C}^p,
\end{eqnarray}
of $p$ linear differential equations with matrix $B(z)$
meromorphic on the Riemann sphere $\overline{\mathbb C}$ and
holomorphic outside the set of singular points $a_1,\ldots,a_n$.

By the {\it monodromy representation} or the {\it monodromy} of
this system we mean the representation
\begin{eqnarray}\label{monodr}
\chi: \pi_1(\overline{\mathbb C}\setminus\{a_1,\ldots,a_n\})\to
GL(p, {\mathbb C})
\end{eqnarray}
of the fundamental group of the punctured sphere in the space of
non-singular complex matrices of size $p$. This representation is
defined as follows. In a neighbourhood of a non-singular point
$z_0$ we consider a {\it fundamental matrix} $Y(z)$ the columns of
which form a basis in the solution space of the system. Analytic
continuation of the matrix $Y(z)$ along an arbitrary loop $\gamma$
outgoing from $z_0$ and lying in $\overline{\mathbb
C}\setminus\{a_1,\ldots,a_n\}$ transforms this matrix into an (in
general different) fundamental matrix $\widetilde Y(z)$. The two
bases are related by means of a non-singular transition matrix
$G_{\gamma}$ corresponding to the loop $\gamma$:
$$
Y(z)=\widetilde Y(z)G_{\gamma}.
$$
The map $\gamma\mapsto G_{\gamma}$ defines the representation
$\chi$ with respect to $Y(z)$. Since the fundamental group of the
punctured sphere is generated by homotopy classes of all simple
loops $\gamma_i$ (each $\gamma_i$ encircles the only singular
point $a_i$, and by convention we assume the loop
$\gamma_1\ldots\gamma_n$ is contractible), the representation
$\chi$ is defined by local {\it monodromy matrices} $G_i$
corresponding to these loops.

A singular point $a_i$ of the system (\ref{syst}) is said to be
{\it Fuchsian} if the matrix differential 1-form $B(z)dz$ has a
simple pole at this point. By Sauvage's theorem (see \cite{Ha},
Th.\,11.1) a Fuchsian singularity is always {\it regular} (i.~e.,
each solution has at most power growth near it), although a
regular singularity is not necessarily Fuchsian. The system
(\ref{syst}) is said to be {\it Fuchsian} if all its singular
points are Fuchsian.

The classical Riemann-Hilbert problem asks for conditions under
which it is possible to construct a Fuchsian system (\ref{syst})
with prescribed singular points $a_1,\ldots,a_n$ and prescribed
monodromy (\ref{monodr}) (in general case the problem has a
negative solution, the counterexample was found by A.\,Bolibrukh).
One knows various sufficient conditions for the affirmative
solution of this problem. One such condition is the irreducibility
of the representation (\ref{monodr}). And by Plemelj's theorem the
problem has always a solution if one allows the point $a_1$ to be
regular rather than Fuchsian (one can learn all this from
\cite{AB}).

Before formulating the generalized Riemann-Hilbert problem (the
GRH-problem) we recall the notions of local holomorphic and
meromorphic transformations, the Poincar\'e rank and the minimal
Poincar\'e rank of the system (\ref{syst}).

If the coefficient matrix $B(z)$ of the system (\ref{syst}) has
the Laurent expansion of the form
$$
B(z)=\frac{B_{-r-1}}{(z-a)^{r+1}}+\ldots+\frac{B_{-1}}{z-a}+B_0+\ldots
\qquad (B_{-r-1}\ne0)
$$
in a neighbourhood of a singularity $a=a_i$ then we will refer to
the integer $r$ as the {\it Poincar\'e rank} of the system at this
point.

A local linear transformation (in a neighbourhood $O_i$ of a point
$a_i$)
$$
y'=\Gamma(z)y
$$
is said to be {\it holomorphic} (or {\it holomorphically
invertible}) if the matrix $\Gamma(z)$ is holomorphic in $O_i$ and
$\det\Gamma(a_i)\ne0$. And this transformation is said to be {\it
meromorphic} (or {\it meromorphically invertible}) if the matrix
$\Gamma(z)$ is meromorphic at $a_i$, holomorphic in
$O_i\setminus\{a_i\}$ and $\det\Gamma(z)\not\equiv0$.

Such transformations take (\ref{syst}) to the system
\begin{eqnarray}\label{syst'}
\frac{dy'}{dz}=B'(z)y', \qquad
B'(z)=\frac{d\Gamma}{dz}\Gamma^{-1}+\Gamma B(z)\Gamma^{-1}.
\end{eqnarray}
Then the systems (\ref{syst}) and (\ref{syst'}) are called {\it
holomorphically (meromorphically) equivalent}.

A holomorphic transformation does not change the Poincar\'e rank
of the original system, while a meromorphic one can increase or
decrease the Poincar\'e rank. The {\it minimal Poincar\'e rank} of
the system (\ref{syst}) at the point $a_i$ is the smallest
Poincar\'e rank of local systems (\ref{syst'}) in the meromorphic
equivalence class of (\ref{syst}) at the point $a_i$.

Now the GRH-problem can be formulated as follows.
\\

{\it Let for each $i=1,\ldots,n$ a local system
\begin{eqnarray}\label{systi}
\frac{dy}{dz}=B_i(z)y
\end{eqnarray}
be given in the neighbourhood $O_i$ of the (irregular) singular
point $a_i$ of the minimal Poincar\'e rank $r_i$, such that its
monodromy matrix coincides with $G_i$. Does there exist a global
system (\ref{syst}) with singularities $a_1,\ldots,a_n$ of the
Poincar\'e ranks $r_1,\ldots,r_n$, with prescribed monodromy
(\ref{monodr}) and such that it is meromorphically equivalent to
the system (\ref{systi}) in each $O_i$?}
\\

Note that the classical Riemann-Hilbert problem also can be
formulated in a such way, but in this case a matrix $B_i(z)$ of a
system (\ref{systi}) in $O_i$ always can be chosen of the form
$B_i(z)=E_i/(z-a_i)$, $E_i=(1/2\pi\sqrt{-1})\ln G_i$. Hence the
systems (\ref{systi}) are uniquely determined by the monodromy
representation (\ref{monodr}) and can be omitted.

We will refer to the monodromy representation (\ref{monodr}) and
local systems (\ref{systi}) as the {\it generalized monodromy
data}.

These data are called {\it reducible} if the representation
(\ref{monodr}) is reducible and the local systems (\ref{systi})
are also reducible, i.~e., they can be reduced via meromorphic
transformations to systems with coefficient matrices of the same
block upper-triangular form. Otherwise we say that the generalized
monodromy data are {\it irreducible}.

A.\,Bolibrukh has generalized his method of solution of the
classical Riemann-Hilbert problem to the case of irregular
singularities (the GRH-problem) and has obtained together with the
co-authors in \cite{BMM} some sufficient conditions for the
affirmative solution of the problem. One such condition is the
irreducibility of the generalized monodromy data in the case if
one at least of the singularities is {\it unramified} (the
definition of ramified and unramified singular points see in \S2).

An analogue of Plemelj's theorem is that the problem has always a
solution if one allows the Poincar\'e rank of a global system at
the point $a_1$ not to be minimal. We obtain here an estimate for
the Poincar\'e rank at this point.

{\bf Theorem 1.} {\it Each generalized monodromy data can be
realized by a global system (\ref{syst}) that has the minimal
Poincar\'e ranks at all points but one ($a_1$ for instance), at
which it has the Poincar\'e rank not greater than
$r_1+(p-1)(n+R-1)$, where $R=\sum_{i=1}^nr_i$.}

We also discuss the problem of the construction of a scalar linear
differential equation
$$
\frac{d^py}{dz^p}+b_1(z)\frac{d^{p-1}y}{dz^{p-1}}+\ldots+b_p(z)y=0
$$
with prescribed singular points $a_1,\ldots,a_n$ and generalized
monodromy data. In the construction there necessary arise {\it
apparent singularities} (at which coefficients of an equation are
singular, but solutions are meromorphic, so that a monodromy is
trivial), the number of which we estimate (Theorem 2).

The results of Theorems 1 and 2 may be interpreted as an extension
of Cor.\,1 and Th.\,2 from \cite{VG} to the case of irregular
singularities.
\\
\\

\centerline{\bf\S2. Irregular systems and holomorphic vector
bundles}
\medskip

In this paragraph we recall the main results that we will need
from the theory of irregular singularities and their relations
with vector bundles. Our main reference is the article \cite{BMM}
by Bolibrukh, Malek, Mitschi.

In a neighbourhood of an irregular singularity $a=a_i$ of
Poincar\'e rank $r$ the system (\ref{syst}) has a formal
fundamental matrix $\widehat Y(z)$ of the form
\begin{eqnarray}\label{ffm}
\widehat Y(z)=\widehat F(z)(z-a)^{E}\,Ue^{Q(z)}
\end{eqnarray}
(see \cite{BJL}, Th.\,1), where

$\widehat F(z)$ is a formal (matrix) Laurent series in $z-a$ (in
general divergent) with finite principal part and $\det\widehat
F(z)$ is distinct from the zero series;

$Q(z)$, $E$ and $U$ are block-diagonal matrices with diagonal
blocks $Q^j(z)$, $E^j$ and $U^j$ of the same size, $j=1,\ldots,N$;

the blocks $Q^j(z)$ and $E^j$ too are block-diagonal of the form
$$
Q^j(z)={\rm diag}\left(q_j(t)I_{m_j},q_j(t\zeta_j)I_{m_j},
                       \ldots,q_j(t\zeta_j^{s_j-1})I_{m_j}\right),
$$
where $q_j(t)$ is a polynomial in $t=(z-a)^{-1/s_j}$ with no
constant term, $\zeta_j=e^{-2\pi i/s_j}$ for some integer $s_j$
and $\deg q_j\leqslant rs_j$, $I_{m_j}$ denotes the identity
matrix of size $m_j$;
$$
E^j={\rm diag}\left(\widehat E^{m_j},\widehat E^{m_j}+
\frac1{s_j}I_{m_j},\ldots,\widehat E^{m_j}+
\frac{s_j-1}{s_j}I_{m_j}\right),
$$
where $\widehat E^{m_j}$ is a constant matrix of size $m_j$ in
canonical Jordan form and its eigenvalues $\rho$ satisfy the
condition $0\leqslant{\rm Re}\,\rho<1/s_j$;

the matrix $U^j$ decomposes into blocks $\left[U^j\right]^{kl}$ of
the form
$$
\left[U^j\right]^{kl}=\zeta_j^{-(k-1)(l-1)}I_{m_j}, \qquad
1\leqslant k,l\leqslant s_j,
$$
with respect to the block structure of the matrices $Q^j(z)$ and
$E^j$.

Let the matrix $Q(z)$ be thought of as the (matrix) polynomial in
$1/(z-a)$ of fractional degree $\deg Q$. Then this degree is
called the {\it Katz rank} of a singularity $z=a$.

Since the matrix $Q(z)$ is a meromorphic invariant of the system
(\ref{syst}), it follows from the properties of this matrix that
the Katz rank is not greater than the minimal Poincar\'e rank of a
singularity. Moreover, the minimal Poincar\'e rank is the least
integer greater than or equal to the Katz rank of a singularity.
\\

{\bf Definition 1.} An irregular singularity of the system
(\ref{syst}) is called {\it unramified} (or a {\it singularity
without roots}) if for every block $Q^j(z)$ of the matrix $Q(z)$
from (\ref{ffm}) the corresponding integer $s_j$ is equal to one.

In an opposite case a singularity is called {\it ramified} (or a
{\it singularity with roots}).
\\

Note that in the unramified case the form of the formal
fundamental matrix $\widehat Y(z)$ is simpler:
\begin{eqnarray}\label{uffm}
\widehat Y(z)=\widehat F(z)(z-a)^Ee^{Q(z)},
\end{eqnarray}
where $Q(z)$ decomposes into a direct sum of scalar blocks
$q_j(1/(z-a))I_{m_j}$ with polynomials $q_j$ of degree non greater
than $r$ (and with at least one $q_j$ of degree exactly $r$) and
$E$ is a direct sum of blocks $E^j$ in Jordan normal form with
eigenvalues $\rho$ satisfying $0\leqslant{\rm Re}\,\rho<1$. One
can see that in this case the Poincar\'e rank coincides with the
Katz rank and with the minimal Poincar\'e rank of a singularity.

Now we will describe briefly a method of solution for the
GRH-problem (for details see \cite{BMM}).

From the representation (\ref{monodr}) one constructs over the
punctured Riemann sphere $\overline{\mathbb C}
\setminus\{a_1,\ldots,a_n\}$ a holomorphic vector bundle $F$ of
rank $p$ with a holomorphic connection $\nabla$ having the
prescribed monodromy (\ref{monodr}). This bundle is defined by a
set $\{U_{\alpha}\}$ of sufficiently small discs covering
$\overline{\mathbb C}\setminus\{a_1,\ldots,a_n\}$ and a set
$\{g_{\alpha\beta}\}$ of constant matrices defining a gluing
cocycle. A connection $\nabla$ is defined by a set
$\{\omega_{\alpha}\}$ of matrix differential 1-forms
$\omega_{\alpha}\equiv0$. So in the intersections $U_{\alpha}\cap
U_{\beta}\ne\varnothing$ the gluing conditions
\begin{eqnarray}\label{glue}
\omega_{\alpha}=(dg_{\alpha\beta})g_{\alpha\beta}^{-1}+
g_{\alpha\beta}\omega_{\beta}g_{\alpha\beta}^{-1}
\end{eqnarray}
hold. The connection defines locally the systems
$dy=\omega_{\alpha}y$. A set $\{s_{\alpha}\}$ of solutions to
these systems that satisfy conditions
$s_{\alpha}=g_{\alpha\beta}s_{\beta}$ on $U_{\alpha}\cap
U_{\beta}\ne\varnothing$ defines a {\it horizontal section} for
the connection. The monodromy of the connection (similarly to the
monodromy of the system (\ref{syst})) describes the branching
pattern of horizontal sections after their analytic continuations
along closed paths encircling the points $a_1,\ldots,a_n$.

Further one extends the pair $(F,\nabla)$ to the whole Riemann
sphere by means of the local matrix differential 1-forms
$\omega_i=B_i(z)dz$ of the coefficients of the systems
(\ref{systi}) defined each in the neighbourhood $O_i$ of the point
$a_i$, $i=1,\ldots,n$. This extension has the following coordinate
description. For each $O_i$ consider a fundamental matrix $Y_i(z)$
of the corresponding system (\ref{systi}) and for a nonempty
intersection $O_i\cap U_{\alpha}$ put $g_{i\alpha}(z)=Y_i(z)$ in
this intersection. For any other $U_{\beta}$ that has a nonempty
intersection with $O_i$ define $g_{i\beta}(z)$ as a suitable
analytic continuation of $g_{i\alpha}(z)$ into $O_i\cap U_{\beta}$
such that the set $\{g_{\alpha\beta}, g_{i\alpha}(z)\}$ defines a
cocycle for the covering $\{U_{\alpha}, O_i\}$ of the sphere.
Thus, one gets a vector bundle $F^0$ over the whole Riemann
sphere. Then the set $\{\omega_{\alpha}, \omega_i\}$ will define a
connection $\nabla^0$ on this bundle, because alongside the gluing
conditions (\ref{glue}) for the nonempty intersections
$U_{\alpha}\cap U_{\beta}$ one has
$$
(dg_{i\alpha})g_{i\alpha}^{-1}+g_{i\alpha}\omega_{\alpha}g_{i\alpha}^{-1}
=(dY_i)Y_i^{-1}=\omega_i,
$$
that is a gluing condition for $O_i\cap U_{\alpha}\ne\varnothing$.

The pair $(F^0,\nabla^0)$ is the so-called {\it canonical}
extension of the pair $(F,\nabla)$ in the sense of Deligne.

Now one can construct a family $\cal F$ of extensions of the pair
$(F,\nabla)$ replacing the functions $g_{i\alpha}(z)$ in the
construction of $(F^0,\nabla^0)$ by the functions
\begin{eqnarray}\label{ext}
g'_{i\alpha}(z)=\Gamma_i(z)g_{i\alpha}(z),
\end{eqnarray}
and the forms $\omega_i$ by the forms
\begin{eqnarray}\label{systi2}
\omega'_i=(d\Gamma_i)\Gamma_i^{-1}+\Gamma_i\omega_i\Gamma_i^{-1},
\end{eqnarray}
where $y'=\Gamma_i(z)y$ are all possible meromorphic
transformations of a system (\ref{systi}) not increasing its
Poincar\'e rank $r_i$, $i=1,\ldots,n$ (see (\ref{syst'})). Then
one easily verifies that the gluing conditions
\begin{eqnarray}\label{gluei}
(dg'_{i\alpha})g'^{-1}_{i\alpha}+g'_{i\alpha}\omega_{\alpha}g'^{-1}_{i\alpha}
=\omega'_i
\end{eqnarray}
hold for the nonempty intersections $O_i\cap U_{\alpha}$.

If some bundle $F'$ from the family $\cal F$ is holomorphically
trivial then its connection $\nabla'$ defines a global system
(\ref{syst}) solving the GRH-problem.

Indeed, the triviality of the bundle $F'$ means that for the
covering $\{U_{\alpha}, O_i\}$ of the sphere there exists a
corresponding set $\{h_{\alpha}(z), h_i(z)\}$ of holomorphically
invertible matrix functions such that
\begin{eqnarray}\label{triv}
h_{\alpha}(z)g_{\alpha\beta}=h_{\beta}(z), \qquad
h_i(z)g'_{i\alpha}(z)=h_{\alpha}(z)
\end{eqnarray}
respectively in $U_{\alpha}\cap U_{\beta}\ne\varnothing$, $O_i\cap
U_{\alpha}\ne\varnothing$. These relations alongside the gluing
conditions (\ref{gluei}) imply that the forms
\begin{eqnarray}\label{globalform}
\tilde\omega_i=(dh_i)h_i^{-1}+h_i\omega'_ih_i^{-1}, \qquad
\tilde\omega_{\alpha}=(dh_{\alpha})h_{\alpha}^{-1}+
h_{\alpha}\omega_{\alpha}h_{\alpha}^{-1}
\end{eqnarray}
coincide over corresponding nonempty intersections and thus define
a global form $\omega=B(z)dz$ on the whole Riemann sphere.

By the construction the global system $dy=\omega y$ is
meromorphically equivalent to the original local systems
(\ref{systi}) in each $O_i$, has the Poincar\'e ranks
$r_1,\ldots,r_n$ and the prescribed monodromy (\ref{monodr}).

From the other side, the existence of the global system solving
the GRH-problem for the generalized monodromy data (\ref{monodr}),
(\ref{systi}) implies the triviality of some bundle from $\cal F$.
\\

Thus, one gets that {\it the GRH-problem has a positive solution
for the given generalized monodromy data (\ref{monodr}),
(\ref{systi}) if and only if at least one of the vector bundles of
the family $\cal F$ is holomorphically trivial.}
\\

The Birkhoff-Grothendieck theorem states that each holomorphic
vector bundle $F'$ of rank $p$ over the Riemann sphere is
holomorphically equivalent to a sum of line bundles
$$
F'\cong{\cal O}(k_1)\oplus\ldots\oplus{\cal O}(k_p),
$$
where $\{k_1\geqslant\ldots\geqslant k_p\}$ is a system of
integers called the {\it splitting type} of the bundle $F'$.

This theorem means that for a cocycle $\{g_{\alpha\beta},
g'_{i\alpha}(z)\}$ defining the bundle $F'$ relations (\ref{triv})
hold, where all matrix functions $h_i(z)$ except one ($h_1(z)$ for
instance) are holomorphically invertible in corresponding $O_i$
and $h_1(z)$ is of the form $h_1(z)=(z-a_1)^{-K}h'_1(z)$, where
$K={\rm diag}(k_1,\ldots,k_p)$ and $h'_1(z)$ is holomorphically
invertible in $O_1$. Thus, one gets a global system $dy=\omega y$
with given generalized monodromy data, but the Poincar\'e rank of
this system at the singular point $a_1$ may be greater than $r_1$,
since
\begin{eqnarray}\label{formsO}
\omega=-\frac K{z-a_1}dz+(z-a_1)^{-K}\tilde\omega_1(z-a_1)^K
\end{eqnarray}
in $O_1$, where ${\rm ord}_{a_1}\tilde\omega_1=-(r_1+1)$. The
relation (\ref{formsO}) implies that the Poincar\'e rank of the
global system at the singularity $a_1$ is not greater than
$r_1+k_1-k_p$. Further to prove Theorem 1 we will estimate
integers $k_1-k_p$ for some bundles $F'$ from $\cal F$.
\\

{\bf Definition 2.} The {\it degree} $\deg F'$ of a bundle $F'$
with a connection $\nabla'$ is the sum
$$
\deg F'=\sum_{i=1}^n{\rm res}_{a_i}{\rm tr}\,\omega'_i
$$
determined by the forms $\omega'_i$ from (\ref{systi2}).
\\

The degree of a bundle is an integer equal to the sum of the
coefficients $k_i$ of its splitting type. Indeed, from the
relations (\ref{globalform}) and (\ref{formsO}) it follows that
\begin{eqnarray*}
{\rm res}_{a_1}{\rm tr}\,\omega &=& -{\rm tr}\,K+{\rm res}_{a_1}
{\rm tr}\,\tilde\omega_1=-{\rm tr}\,K+{\rm res}_{a_1}{\rm tr}\,
\omega'_1, \\
{\rm res}_{a_i}{\rm tr}\,\omega &=& {\rm res}_{a_i}{\rm
tr}\,\tilde\omega_i={\rm res}_{a_i}{\rm tr}\,\omega'_i, \qquad
i=2,\ldots,n.
\end{eqnarray*}
Then by theorem on the sum of residues one has $\sum_{i=1}^n{\rm
res}_{a_i}{\rm tr}\,\omega=0$ and therefore
$$
\deg F'=\sum_{i=1}^n{\rm res}_{a_i}{\rm tr}\,\omega'_i={\rm
tr}\,K.
$$

Now we consider a subset ${\cal E}\subset\cal F$ of the family
$\cal F$ constructed by means of meromorphic transformations with
matrices $\Gamma_i(z)$ from (\ref{ext}), (\ref{systi2}) of some
special form. For this construction one needs the following
definition.
\\

{\bf Definition 3.} Consider a system (\ref{systi}) with an
(irregular) singular point $a_i$ and its formal fundamental matrix
$\widehat Y_i(z)$ of the form (\ref{ffm}), where all matrices are
supplied with subscript $i$. An {\it admissible matrix} for this
system is an integer-valued diagonal matrix $\Lambda_i={\rm
diag}(\Lambda_i^1,\ldots,\Lambda_i^N)$ blocked in the same way as
$Q_i(z)$ and such that

$(z-a_i)^{\Lambda_i^j}E_i^j(z-a_i)^{-\Lambda_i^j}$ is holomorphic
at the point $a_i$ if the block $Q_i^j(z)$ has no ramification;

$\Lambda_i^j$ is a scalar matrix if the block $Q_i^j(z)$ has
ramification.
\\

{\bf Remark 1.} If a singularity $a_i$ is Fuchsian for a system
(\ref{systi}), i.~e., $a_i$ is a regular singularity with minimal
Poincar\'e rank, then there exists a fundamental matrix $Y_i(z)$
of this system of the form
\begin{eqnarray}\label{fmi}
Y_i(z)=M(z)(z-a_i)^{E_i},
\end{eqnarray}
where the matrix $M(z)$ is meromorphic at the point $a_i$ and
$E_i$ is an upper-triangular constant matrix with eigenvalues
$\rho$ satisfying $0\leqslant{\rm Re}\,\rho<1$.

In this case an admissible matrix $\Lambda_i$ is a diagonal
integer-valued matrix such that the matrix
$(z-a_i)^{\Lambda_i}E_i(z-a_i)^{-\Lambda_i}$ is holomorphic at the
point $a_i$.
\\

Let us write the matrix $\widehat Y_i(z)$ as follows:
\begin{eqnarray}\label{ffmi}
\widehat Y_i(z)=\widehat
F_i(z)(z-a_i)^{-\Lambda_i}(z-a_i)^{\Lambda_i}(z-a_i)^{E_i}\,U_ie^{Q_i(z)}.
\end{eqnarray}
By analogue of Sauvage's lemma (see \cite{Ha}, L.\,11.2) for
formal series, there exists a meromorphically invertible matrix
$\Gamma'_i(z)$ in $O_i$, such that
\begin{eqnarray}\label{sauvage}
\Gamma'_i(z)\widehat F_i(z)(z-a_i)^{-\Lambda_i}=(z-a_i)^D\widehat
F_0(z),
\end{eqnarray}
where $D$ is a diagonal integer-valued matrix and $\widehat
F_0(z)$ is an invertible formal (matrix) Taylor series in $z-a_i$.

The required meromorphic transformations for an irregular
singularity $a_i$ are now defined by the matrices
$\Gamma^{\Lambda_i}(z)= (z-a_i)^{-D}\,\Gamma'_i(z)$ depending on
$\Lambda_i$ (because $\Gamma'_i(z)$ depends on $\Lambda_i$), and
by the matrices $(z-a_i)^{\Lambda_i}M^{-1}(z)$ for a Fuchsian
singularity $a_i$, $i=1,\ldots,n$.

One needs only to verify that such transformations do not increase
the Poincar\'e ranks $r_i$ of the systems (\ref{systi}). This is
provided by the following lemma.

{\bf Lemma 1.} {\it Consider the set $\cal E$ of the extensions
$(F^{\Lambda},\nabla^{\Lambda})$ of the pair $(F,\nabla)$ to the
whole Riemann sphere obtained by means of all possible systems
$\Lambda=\{\Lambda_1,\ldots,\Lambda_n\}$ of admissible matrices
for the singularities $a_1,\ldots,a_n$. This set is a subset of
the family $\cal F$, i.~e., for each pair
$(F^{\Lambda},\nabla^{\Lambda})$ the Poincar\'e rank of the
connection $\nabla^{\Lambda}$ at the point $a_i$ is equal to
$r_i$.

Moreover, for the degree of the bundle $F^{\Lambda}$ the following
relation holds:
$$
\deg F^{\Lambda}=\sum_{i=1}^n{\rm tr}\,(\Lambda_i+E_i).
$$
}

{\bf Proof.} In the case of irregular singularity $a_i$ it follows
from (\ref{ffmi}), (\ref{sauvage}) that the transformation
$y'=\Gamma^{\Lambda_i}(z)y$ transforms the system (\ref{systi}) to
the system with formal fundamental matrix $\widehat Y'_i(z)$ of
the form
$$
\widehat Y'_i(z)=\Gamma^{\Lambda_i}(z)\widehat Y_i(z)=\widehat
F_0(z)(z-a_i)^{\Lambda_i}(z-a_i)^{E_i}\,U_ie^{Q_i(z)},
$$
therefore its Poincar\'e rank is equal to the Poincar\'e rank of
the system with fundamental matrix $\widetilde
Y(z)=(z-a_i)^{\Lambda_i}(z-a_i)^{E_i}\,U_ie^{Q_i(z)}$ (because the
formal holomorphic transformation $\tilde y=\widehat
F_0^{-1}(z)y'$ does not change the Poincar\'e rank) and with
coefficient matrix
\begin{eqnarray*}
\widetilde B(z) & = &\frac{d\widetilde Y}{dz}\widetilde Y^{-1} =
\frac{\Lambda_i}{z-a_i}+(z-a_i)^{\Lambda_i}\frac {E_i}{z-a_i}(z-a_i)^{-\Lambda_i}+\\
& & +(z-a_i)^{\Lambda_i}(z-a_i)^{E_i}\,U_i\,\frac{dQ_i}{dz}\,
U_i^{-1}(z-a_i)^{-E_i}(z-a_i)^{-\Lambda_i}.
\end{eqnarray*}
This matrix is blocked in the same way as the matrices $Q_i(z)$,
$E_i$ and $U_i$.

If a block $Q_i^j(z)$ has no ramification then it is scalar,
$(z-a_i)^{\Lambda_i^j}{E_i^j}(z-a_i)^{-\Lambda_i^j}$ is
holomorphic at the point $a_i$ and $U_i^j=I$, therefore
$$
\widetilde
B^j(z)=\frac{\Lambda_i^j}{z-a_i}+(z-a_i)^{\Lambda_i^j}\frac
{E_i^j}{z-a_i}(z-a_i)^{-\Lambda_i^j}+\frac{dQ_i^j}{dz}.
$$

If a block $Q_i^j(z)$ has ramification then the matrix
$\Lambda_i^j$ is scalar, therefore
$$
\widetilde B^j(z)=\frac{\Lambda_i^j}{z-a_i}+\frac{E_i^j}{z-a_i}
+(z-a_i)^{E_i^j}\,U_i^j\,\frac{dQ_i^j}{dz}\,(U_i^j)^{-1}(z-a_i)^{-E_i^j}.
$$

The (generally speaking, fractional) degree of the polynomial
$\frac{dQ_i}{dz}$ in $1/(z-a_i)$ is at most $r_i+1$, and the real
parts of the eigenvalues of the matrix $E_i$ lie in the half-open
interval $[0, 1)$, therefore in any case ${\rm
ord}_{a_i}\widetilde B(z)=-(r_i+1)$ (it follows from the fact that
$\widetilde B(z)$ contains only integer powers of $z-a_i$).

In the case of Fuchsian singularity $a_i$ it follows from
(\ref{fmi}) that the transformation
$y'=(z-a_i)^{\Lambda_i}M^{-1}(z)y$ takes the system (\ref{systi})
to the system with fundamental matrix $Y'_i(z)$ of the form
$$
Y'_i(z)=(z-a_i)^{\Lambda_i}M^{-1}(z)Y_i(z)=(z-a_i)^{\Lambda_i}(z-a_i)^{E_i},
$$
therefore the new system is also Fuchsian at the point $a_i$.

Now let us prove the second part of the lemma. By definition $\deg
F^{\Lambda}=\sum_{i=1}^n{\rm res}_{a_i}{\rm
tr}\,\omega^{\Lambda_i}$, where
\begin{eqnarray*}
\omega^{\Lambda_i}&=&(d\widehat Y'_i)\widehat Y'^{-1}_i=(d\widehat
F_0)\widehat F_0^{-1}+ \\
& &+\widehat F_0\left(\frac{\Lambda_i}{z-a_i}+
(z-a_i)^{\Lambda_i}\frac {E_i}{z-a_i}(z-a_i)^{-\Lambda_i}\right.+\\
& &
+\left.(z-a_i)^{\Lambda_i}(z-a_i)^{E_i}\,U_i\,\frac{dQ_i}{dz}\,
U_i^{-1}(z-a_i)^{-E_i}(z-a_i)^{-\Lambda_i}\right)\widehat
F_0^{-1}dz
\end{eqnarray*}
if $a_i$ is irregular, and
\begin{eqnarray*}
\omega^{\Lambda_i}=(dY'_i)Y'^{-1}_i=\frac{\Lambda_i}{z-a_i}dz+
(z-a_i)^{\Lambda_i}\frac {E_i}{z-a_i}(z-a_i)^{-\Lambda_i}dz
\end{eqnarray*}
if $a_i$ is Fuchsian.

Thus, one can see that in both cases ${\rm res}_{a_i}{\rm
tr}\,\omega^{\Lambda_i}={\rm tr}\,(\Lambda_i+E_i)$ and therefore
$\deg F^{\Lambda}=\sum_{i=1}^n{\rm tr}\,(\Lambda_i+E_i)$. {\hfill
$\Box$}

Let us call the eigenvalues $\beta_i^j=\lambda_i^j+\rho_i^j$ of
the matrix $\Lambda_i+E_i$ {\it (formal) exponents of the
connection $\nabla^{\Lambda}$} at the (irregular) singular point
$a_i$.
\\
\\

\centerline{\bf\S3. Proof of Theorem 1}
\medskip

Theorem 1 is a direct consequence of the following result (which
is based on the proof of L.\,2 from \cite{BMM}).

{\bf Proposition 1.} {\it Consider a pair
$(F^{\Lambda},\nabla^{\Lambda})\in\cal E$ such that the exponents
of $\nabla^{\Lambda}$ satisfy the condition $0\leqslant{\rm
Re}\,\beta_i^j<M$, $M\in\mathbb N$.

Then the following inequalities hold for the splitting type
$(k_1^{\Lambda},\ldots,k_p^{\Lambda})$ of the bundle
$F^{\Lambda}$:
$$
k_j^{\Lambda}-k_{j+1}^{\Lambda}\leqslant(n+R)M-1,\qquad
j=1,\ldots,p-1,
$$
where $R=\sum_{i=1}^nr_i$.}

{\bf Proof.} We consider two separate cases.

{\it Case 1.} For the splitting type of the bundle $F^{\Lambda}$
one has the inequalities
$$
k_j^{\Lambda}-k_{j+1}^{\Lambda}\leqslant n+R-2,\qquad
j=1,\ldots,p-1.
$$
Since $M\in\mathbb N$, the required result in this case follows
from these inequalities.

{\it Case 2.} For some $l$ one has
$k_l^{\Lambda}-k_{l+1}^{\Lambda}>n+R-2$.

Consider a system (\ref{syst}) with singularities $a_1,\ldots,a_n$
and generalized monodromy data (\ref{monodr}), (\ref{systi}) such
that the Poincar\'e ranks of singularities $a_2,\ldots,a_n$ are
equal to $r_2,\ldots,r_n$ respectively and the differential 1-form
$\omega=B(z)dz$ of the coefficients in the neighbourhood $O_1$ of
the point $a_1$ has the form
\begin{eqnarray}\label{formsO2}
\omega=-\frac
K{z-a_1}dz+(z-a_1)^{-K}\tilde\omega^{\Lambda_1}(z-a_1)^K,
\end{eqnarray}
where $K={\rm diag}(k_1^{\Lambda},\ldots,k_p^{\Lambda})$ and ${\rm
ord}_{a_1}\tilde\omega^{\Lambda_1}=-(r_1+1)$ (see (\ref{formsO})).

By (\ref{formsO2}) the entries $\omega_{mj}$ and
$\tilde\omega_{mj}$ of the matrix differential 1-forms $\omega$
and $\tilde\omega^{\Lambda_1}$ are connected for $m\ne j$ by the
equality
$$
\omega_{mj}=(z-a_1)^{-k_m^{\Lambda}+k_j^{\Lambda}}\,\tilde\omega_{mj}.
$$
By assumption $k_l^{\Lambda}-k_{l+1}^{\Lambda}>n+R-2$ for some
$l$, therefore we have $k_j^{\Lambda}-k_m^{\Lambda}>n+R-2$ for
$j\leqslant l$, $m>l$. Hence the orders ${\rm ord}_{a_1}
\omega_{mj}$ at the point $a_1$ of the differential 1-forms
$\omega_{mj}$ with indicated indices are greater than $n+R-r_1-3$,
whereas the sum of the orders ${\rm ord}_{a_i}\omega_{mj}$ at the
singular points distinct from $a_1$ is at least $-n-R+r_1+1$.

We thus obtain for meromorphic forms $\omega_{mj}$ with indicated
indices that the sum of their orders over all singularities and
zeros is greater than $-2$, although this sum is known to be $-2$
for a non-trivial differential 1-form on $\overline{\mathbb C}$
(the degree of the canonical divisor; see \cite{Fo},
Prop.\,17.12). Hence these forms are identically equal to zero, so
that in an equivalent coordinate description of the bundle
$F^{\Lambda}$ (defined by a cocycle $\{\tilde g_{\alpha\beta},
\tilde g_{i\alpha}\}$ with only functions of the form $\tilde
g_{1\alpha}(z)=(z-a_1)^K$ distinct from the identity matrix) all
matrix differential 1-forms $\tilde\omega^{\Lambda_i}$,
$\tilde\omega_{\alpha}$ defining the connection $\nabla^{\Lambda}$
are block upper-triangular:
$$
\tilde\omega^{\Lambda_i}=\begin{pmatrix} \tilde\omega_i^1 & *\\
                                        0 & \tilde\omega_i^2
                         \end{pmatrix}, \qquad
\tilde\omega_{\alpha}=   \begin{pmatrix} \tilde\omega_{\alpha}^1 & *\\
                                        0 & \tilde\omega_{\alpha}^2
                         \end{pmatrix},
$$
where all matrix forms $\tilde\omega_i^1$,
$\tilde\omega_{\alpha}^1$ have size $l\times l$.

This means that the bundle $F^{\Lambda}$ has a subbundle
$F^1\cong{\cal O}(k_1^{\Lambda})\oplus\ldots\oplus{\cal
O}(k_l^{\Lambda})$ of rank $l$ with a connection $\nabla^1$
defined by the forms $\tilde\omega_i^1$, $\tilde\omega_{\alpha}^1$
satisfying the required gluing conditions (in view of
(\ref{formsO2})).

From results of \cite{BJL2} it follows that a formal fundamental
matrix $\widetilde Y_i$ of a local system
$dy=\tilde\omega^{\Lambda_i}y$ (which is holomorphically
equivalent to a system $dy=\omega^{\Lambda_i}y$) with irregular
singularity $a_i$ can be chosen to have a block upper-triangular
structure similar to $\tilde\omega^{\Lambda_i}$:
$$
\widetilde Y_i=\begin{pmatrix} \widetilde Y_i^1 & *\\
                               0 & \widetilde Y_i^2
               \end{pmatrix},
$$
furthermore it has the form
$$
\widetilde Y_i(z)=\widetilde F_0(z)(z-a_i)^{\widetilde\Lambda_i}
(z-a_i)^{\widetilde E_i}\,\widetilde U_ie^{\widetilde Q_i(z)},
$$
where $\widetilde\Lambda_i=S^{-1}\Lambda_iS$, $\widetilde
E_i=S^{-1}E_iS$, $\widetilde U_i=S^{-1}U_iS$, $\widetilde
Q_i(z)=S^{-1}Q_i(z)S$ for some constant invertible matrix $S$, the
matrices $\widetilde\Lambda_i$ and $\widetilde Q_i(z)$ are
diagonal and obtained by suitable permutations of the diagonal
elements of $\Lambda_i$ and $Q_i(z)$ respectively, the matrix
$\widetilde E_i$ is upper-triangular and the matrix $\widetilde
U_i={\rm diag}(\widetilde U_i^1, \widetilde U_i^2)$ is block
diagonal with respect to the block structure of the matrix
$\widetilde Y_i$. Moreover, the invertible formal (matrix) Taylor
series $\widetilde F_0(z)$ has the same block upper-triangular
structure as the matrix $\widetilde Y_i$.

Thus, one gets that the set $\{{^1\beta_i^1},\ldots,{^1\beta_i^l}
\}$ of the (formal) exponents of the connection $\nabla^1$ at the
singularity $a_i$ is a subset of the (formal) exponents of the
connection $\nabla^{\Lambda}$ at this point.

Assume that $k_l^{\Lambda}-k_{l+1}^{\Lambda}\geqslant(n+R)M$. Then
for the mean value of the exponents $^1\beta_i^j$ of the
connection $\nabla^1$ we have the lower bound
$$
\frac1{ln}\sum_{i=1}^n\sum_{j=1}^l{^1\beta_i^j}=\frac{\deg
F^1}{ln}=\frac{k_1^{\Lambda}+\ldots+k_l^{\Lambda}}{ln}
\geqslant\frac{k_{l+1}^{\Lambda}}n+M,
$$
while for the mean value of the other exponents $^2\beta_i^j$ of
the connection $\nabla^{\Lambda}$ we have the upper bound
$$
\frac1{(p-l)n}\sum_{i=1}^n\sum_{j=1}^{p-l}{^2\beta_i^j}=\frac{\deg
F^{\Lambda}-\deg F^1}{(p-l)n}=\frac{k_{l+1}^{\Lambda}+\ldots+
k_p^{\Lambda}}{(p-l)n}\leqslant\frac{k_{l+1}^{\Lambda}}n.
$$

Hence the mean value of the exponents $^1\beta_i^j$ is larger by
$M$ at least than the mean value of the exponents $^2\beta_i^j$,
while by the hypothesis the real parts of all the exponents of the
connection $\nabla^{\Lambda}$ are strictly less than $M$. We
arrive to a contradiction, therefore
$k_l^{\Lambda}-k_{l+1}^{\Lambda}\leqslant(n+R)M-1$ for each $l$.
{\hfill $\Box$}
\\

{\bf Proof of Theorem 1.} Consider the pair
$(F^{\Lambda^0},\nabla^{\Lambda^0})\in{\cal E}\subset\cal F$
corresponding to the system $\Lambda^0=\{0,\ldots,0\}$ of zero
matrices. In that case the exponents $\beta_i^j$ of the connection
$\nabla^{\Lambda^0}$ satisfy the condition
$$
0\leqslant{\rm Re}\,\beta_i^j={\rm Re}\,\rho_i^j<1,
$$
therefore by Proposition 1 we have the inequalities
$$
k_j^0-k_{j+1}^0\leqslant n+R-1, \qquad j=1,\ldots,p-1,
$$
for the coefficients $k^0_j$ of the splitting type of the bundle
$F^{\Lambda^0}$. Hence
$$
k_1^0-k_p^0=\sum_{j=1}^{p-1}(k_j^0-k_{j+1}^0)\leqslant
(p-1)(n+R-1).
$$

The coefficient matrix $B(z)$ of the global system (\ref{syst})
corresponding to the connection $\nabla^{\Lambda^0}$ has the form
$$
B(z)=-\frac{K^0}{z-a_1}+(z-a_1)^{-K^0}\widetilde B(z)(z-a_1)^{K^0}
$$
in the neighbourhood $O_1$ of the point $a_1$, where $K^0={\rm
diag}(k_1^0,\ldots,k_p^0)$ and ${\rm ord}_{a_1}\widetilde
B(z)=-(r_1+1)$ (see (\ref{formsO2})). Then the Poincar\'e rank of
this system at the point $a_1$ is not greater than the quantity
$$
r_1+k_1^0-k_p^0\leqslant r_1+(p-1)(n+R-1)
$$
(recall that this system has the prescribed singularities
$a_1,\ldots,a_n$, generalized monodromy data (\ref{monodr}),
(\ref{systi}) and the Poincar\'e ranks $r_2,\ldots,r_n$ at the
points $a_2,\ldots,a_n$ respectively). {\hfill $\Box$}
\\

{\bf Remark 2.} From the proof of Proposition 1 it follows that if
one at least of the singularities $a_1,\ldots,a_n$ is irregular
(and then $R>0$) then
$k_j^{\Lambda}-k_{j+1}^{\Lambda}\leqslant(n+R)M-2$ for each $j$.

In this case one gets that in Theorem 1 the Poincar\'e rank of the
global system at the point $a_1$ is not greater than
$r_1+(p-1)(n+R-2)$.
\\

Let us tell some words about the problem of the meromorphic
transformation of a system
\begin{eqnarray}\label{birksyst}
\frac{dy}{dz}=C(z)y, \qquad C(z)=\frac{C_{-r-1}}{z^{r+1}}+
\ldots+\frac{C_{-1}}z+C_0+\ldots,
\end{eqnarray}
of $p$ linear differential equations to a {\it Birkhoff standard
form} in a neighbourhood of an irregular singularity $z=0$ of
Poincar\'e rank $r$ (not necessarily minimal), i.e., to a system
with coefficient matrix $C'(z)$ of the form
\begin{eqnarray}\label{birk}
C'(z)=\frac{C'_{-r'-1}}{z^{r'+1}}+\ldots+\frac{C'_{-1}}z, \qquad
r'\leqslant r
\end{eqnarray}
(note that such system is defined on the whole Riemann sphere and
$\infty$ is a Fuchsian singularity for it).

This problem is not yet resolved, though it is known that the
problem has an affirmative answer in dimensions $p=2$ and $p=3$;
one also knows various sufficient conditions for a positive
solution in an arbitrary dimension $p$ (for instance, the problem
has a positive solution if the system (\ref{birksyst}) is
irreducible (A.\,Bolibrukh) or if all the eigenvalues of the
matrix $C_{-r-1}$ are distinct (H.\,L.\,Turrittin); see Balser's
survey \cite{Ba} for details).

Denote by $r^0>0$ the minimal Poincar\'e rank of the system
(\ref{birksyst}) and consider the GRH-problem for the following
generalized monodromy data:

i) an irregular singularity $a_1=0$ with local system
meromorphically equivalent to (\ref{birksyst}) and of Poincar\'e
rank $r^0$;

ii) a Fuchsian singularity $a_2=\infty$.

By Theorem 1 (where $n=2$ and $R=r^0>0$) and Remark 2 there exists
a global system on the whole Riemann sphere that is Fuchsian at
infinity and meromorphically equivalent to the system
(\ref{birksyst}) in a neighbourhood of the point $a_1=0$. The
coefficient matrix of this system has the form (\ref{birk}), where
$r'\leqslant r^0+(p-1)r^0= pr^0$. Thus, one gets the following
statement.

{\bf Corollary 1.} {\it If for the minimal Poincar\'e rank $r^0$
of the system (\ref{birksyst}) the inequality $r^0\leqslant r/p$
holds then it can be meromorphically transformed to a Birkhoff
standard form.}
\\
\\

\centerline{\bf\S4. The GRH-problem for scalar linear differential
equations}
\medskip

Consider a linear differential equation
\begin{eqnarray}\label{eq}
\frac{d^pu}{dz^p}+b_1(z)\frac{d^{p-1}u}{dz^{p-1}}+\ldots+b_p(z)u=0
\end{eqnarray}
of order $p$ with coefficients $b_1(z),\ldots,b_p(z)$ meromorphic
on the Riemann sphere $\overline{\mathbb C}$ and holomorphic
outside the set of singular points $a_1,\ldots,a_n$.

One defines the monodromy representation
\begin{eqnarray}\label{monodreq}
\chi: \pi_1(\overline{\mathbb C}\setminus\{a_1,\ldots,a_n\})\to
GL(p, {\mathbb C})
\end{eqnarray}
of this equation in the same way as for a system (\ref{syst}); one
merely needs to consider in place of a fundamental matrix $Y(z)$ a
row $(u_1,\ldots,u_p)$, where the functions $u_1(z),\ldots,u_p(z)$
form a basis in the solution space of the equation. This
representation is defined by local monodromy matrices $G_i$
corresponding to simple loops $\gamma_i$.

A singular point $a_i$ of the equation (\ref{eq}) is said to be
{\it Fuchsian} if the coefficient $b_j(z)$ has at this point a
pole of order $j$ or lower ($j=1,\ldots,p$). By Fuchs's theorem
(see \cite{Ha}, Th.\,12.1) a singular point of the equation
(\ref{eq}) is Fuchsian if and only if it is regular. The equation
(\ref{eq}) is said to be {\it Fuchsian} if all its singular points
are Fuchsian.

Using a standard change
$$
y^1=u,\quad y^2=\frac{du}{dz},\quad\ldots,\quad
y^p=\frac{d^{p-1}u}{dz^{p-1}}
$$
one can go over from the equation (\ref{eq}) to a {\it companion
system} (\ref{syst}) with coefficient matrix $B(z)$ of the form
\begin{eqnarray}\label{matreq}
B(z)=\begin{pmatrix} 0 &    1   &        & 0 \\
                       & \ddots & \ddots &   \\
                     0 &        &    0   & 1 \\
                  -b_p & \ldots & \ldots & -b_1 \\
     \end{pmatrix}.
\end{eqnarray}
\\

{\bf Definition 4.} Let us call two linear differential equations
{\it meromorphically equivalent} in a neighbourhood of a singular
point if companion systems are.
\\

The Katz rank $K_i$ of the equation (\ref{eq}) at a singularity
$a_i$ is equal to the Katz rank of the companion system at this
point and it is known (see \cite{LS}, section 3, especially
Th.\,3.2 and Th.\,3.3) that
\begin{eqnarray}\label{katz}
{\rm ord}_{a_i}b_j(z)\geqslant -j(K_i+1),\qquad j=1,\ldots,p.
\end{eqnarray}

Let us now formulate the GRH-problem for scalar linear
differential equations as follows.
\\

{\it Let for each $i=1,\ldots,n$ a local equation
\begin{eqnarray}\label{eqi}
\frac{d^pu}{dz^p}+b^i_1(z)\frac{d^{p-1}u}{dz^{p-1}}+\ldots+b^i_p(z)u=0
\end{eqnarray}
be given in the neighbourhood $O_i$ of the singular point $a_i$ of
Katz rank $K_i$, such that its monodromy matrix coincides with
$G_i$. Does there exist a global equation (\ref{eq}) with
singularities $a_1,\ldots,a_n$, prescribed monodromy
(\ref{monodreq}) and such that it is meromorphically equivalent to
the equation (\ref{eqi}) in each $O_i$?}
\\

We will again refer to the monodromy representation
(\ref{monodreq}) and local equations (\ref{eqi}) as the
generalized monodromy data.

Note that in view of (\ref{katz}) coefficients $b_j(z)$ of a
global equation solving the GRH-problem have bounded orders of
poles.

If $a_i$ is a Fuchsian singularity for the local equation
(\ref{eqi}) for each $i=1,\ldots,n$ then a global equation (if it
exists) is Fuchsian by Fuchs's theorem (or, if the reader prefers,
by (\ref{katz}) and the equalities $K_i=0$). Thus, in this case
one gets a classical problem of the construction of Fuchsian
equation with prescribed singularities and monodromy. Even in this
case the problem has a negative solution in general because for
$p>2$, $n>2$, and for $p=2$, $n>3$ the number of parameters
determining a Fuchsian equation is less than the number of
parameters determining the set of conjugacy classes of
representations $\chi$ (see \cite{AB}, pp. 158--159). Therefore,
to construct such an equation with given monodromy, one needs
so-called apparent singular points. In the case of {\it
irreducible} representation an expression for the smallest
possible number of apparent singular points has been obtained by
A.\,Bolibrukh \cite{Bo}. Some estimate for this number in the case
of arbitrary monodromy is presented in \cite{VG}. Here we extend
this estimate to the case of non-Fuchsian singularities.

{\bf Theorem 2.} {\it Each generalized monodromy data
(\ref{monodreq}), (\ref{eqi}) can be realized by an equation
(\ref{eq}) such that the number of its apparent singularities is
not greater than
$$
\frac{(K+n+1)p(p-1)}2+1,
$$
where $K=-\sum_{i=1}^n[-K_i]$ and $[\;]$ stands for the integer
part.}

{\bf Proof.} From results of J.\,Plemelj it follows that the
classical Riemann-Hilbert problem (for linear systems) has a
positive solution if one at least of the monodromy matrices of the
representation (\ref{monodreq}) is diagonalisable (see \cite{AB},
p. 10, p. 62). Thus, each monodromy representation can be realized
by a Fuchsian system with one apparent singularity (one only needs
to consider the representation $\chi^*$ obtained from
(\ref{monodreq}) by the addition of a singular point $a_{n+1}$
with identity monodromy matrix). In the same way one obtains that
each generalized monodromy data can be realized by a global system
with prescribed singularities of minimal Poincar\'e ranks and
apparent Fuchsian singularity.

Let us consider the representation (\ref{monodreq}) and the local
companion systems for the local equations (\ref{eqi}). For each
local system (with Katz rank $K_i$) consider a meromorphically
equivalent one
\begin{eqnarray}\label{systeqi}
\frac{dy}{dz}=B'_i(z)y
\end{eqnarray}
with minimal Poincar\'e rank $r_i$. Recall that $r_i$ is the least
integer greater than or equal to the Katz rank $K_i$, i.~e.,
$r_i=-[-K_i]$, where $[\;]$ stands for the integer part.

Realize the generalized monodromy data (\ref{monodreq}),
(\ref{systeqi}) by a global system with singularities
$a_1,\ldots,a_n$ of Poincar\'e ranks $r_1,\ldots,r_n$ respectively
and apparent Fuchsian singularity $a_{n+1}$. By Deligne's lemma
(\cite{De}, p. 163) this system can be meromorphically transformed
(globally) to a system with coefficient matrix $B(z)$ of the form
(\ref{matreq}), where $b_1(z),\ldots,b_p(z)$ are meromorphic
functions on the Riemann sphere. Besides $a_1,\ldots,a_{n+1}$, the
transformed system has apparent singularities, the number $m$ of
which satisfies the inequality
$$
m\leqslant\frac{(R+n+1)p(p-1)}2,
$$
where $R=\sum_{i=1}^nr_i$ (this estimate is presented in L.\,2
from \cite{VG}).

One readily sees that the first component of a solution to the
last system is a solution to an equation (\ref{eq}). By the
construction this equation has the prescribed singularities
$a_1,\ldots,a_n$, monodromy (\ref{monodreq}) and it is
meromorphically equivalent to the equation (\ref{eqi}) in each
$O_i$. Furthermore, the number of its apparent singularities is
$m+1$ (note that $a_{n+1}$ is also an apparent singularity of the
equation with respect to the originally prescribed singular points
$a_1,\ldots,a_n$). Bearing in mind that $R=-\sum_{i=1}^n[-K_i]$,
we obtain the required estimate. {\hfill $\Box$}

\end{document}